\documentclass[12pt,amsfonts]{article}
\usepackage[margin=1in]{geometry}  
\usepackage{graphicx}              
\usepackage{amsmath}               
\usepackage{amsfonts}              
\usepackage{amsthm}                

\usepackage{color}

\newcommand{\ud}{\mathrm{d}}
\newcommand{\N}{\mathbb{N}}
\newcommand{\dive}{\mathrm{div}}

\def\R{\mathbb{R}}
\newtheorem{claim}{Claim}[section]
\newtheorem{thm}{Theorem}[section]

\newtheorem{lem}{Lemma}[section]

\newtheorem{cor}{Corollary}[section]
\newcommand{\fim}{\hfill\rule{2mm}{2mm}}

\numberwithin{equation}{section}

\begin{document}
\title{Existence and concentration of solution for a class of fractional  elliptic equation in $\mathbb{R}^N$ via penalization method  }
\author{\sf Claudianor O. Alves\thanks{Research of C. O. Alves partially supported by  CNPq/Brazil  Proc. 304036/2013-7  and INCTMAT/CNPq/Brazi. } 
\\
\small{Universidade Federal de Campina Grande, }\\
\small{Unidade Acad\^emica de Matem\'atica } \\
\small{CEP: 58429-900 - Campina Grande-PB, Brazil}\\
\small{ e-mail: coalves@dme.ufcg.edu}\\
\vspace{1mm}\\
{\sf Ol\'impio H. Miyagaki\thanks{ O.H.M. was partially supported by INCTMAT/CNPq/Brazil, CNPq/Brazil  Proc . 304014/2014-9 and CAPES/Brazil 
Proc. 2531/14-3. 
This paper was completed while the second  author named was visiting the Department of Mathematics of the Rutgers University,
 whose hospitality he gratefully acknowledges. He would like to express his gratitude to Professor Haim Brezis  and Professor Yan Yan Li for invitation and friendship.}} \\
{ \small Departmento de Matem\'atica}\\
{\small Universidade Federal de  Juiz de Fora}\\
{ \small CEP:  36036-330 - Juiz de Fora-MG, Brazil}\\
 {\small e-mail: ohmiyagaki@gmail.com}\vspace{1mm} \\
}
\date{}
\maketitle

\begin{abstract}
In this paper, we study  the existence  and concentration of positive solution for the following class of fractional elliptic equation 
$$
\epsilon^{2s} (-\Delta)^{s}{u}+V(z)u=f(u)\,\,\, \mbox{in} \,\,\, \mathbb{R}^{N},
$$
where  $\epsilon$ is a positive parameter, $f$ has a subcritical growth, $V$ possesses a local minimum,  $N > 2s,$ $s \in (0,1),$ 
and  $ (-\Delta)^{s}u$ is the fractional laplacian. 
\end{abstract}

\vspace{0.5 cm}


\section{Introduction}

The nonlinear fractional  Schr\"{o}dinger
equation
$$
i\epsilon \displaystyle \frac{\partial \Psi}{\partial t}=\epsilon^{2s}(-\Delta)^s
\Psi+(V(z)+E)\Psi-f(\Psi)\,\,\, \mbox{for all}\,\,\, z \in
\R^N,\eqno{(NLS)}
$$
where $N > 2s,$ \ $s\in (0,1),$ \ $\epsilon > 0$, $E\in \R$, and  $V,f$ are continuous functions, has been studied in recent years by many researchers. The
 standing waves solutions of $(NLS)$, namely, $\Psi(z,t)=exp(-i Et)u(z),$ where $u$ is a solution of the
fractional elliptic equation
$$
\ \  \left\{
\begin{array}{l}
\epsilon^{2s} (- \Delta)^su + V(z)u=f(u)
\ \ \mbox{in} \ \ \R^N,
\\
u \in H^s(\R^N),\quad u > 0 \ \ \mbox{on} \ \ \R^N.
\end{array}
\right.
\eqno{(P_{\epsilon})}
$$
In the local case, that is, when $s=1,$ the existence and concentration of
positive solutions for general semilinear elliptic equations
$(P_{\epsilon})$, for the case $N \geq 2$,  have been extensively
studied,  for example,  by  \cite{AOS,Alves10,BPW,BW0,Pino,PFM,OS,FW,Oh1,R,W},  and their references. 
Rabinowitz in \cite{R}  proved the
existence of positive solutions of $(P_{\epsilon})$, for $\epsilon > 0$
small, imposing a global condition,
$$
\liminf_{|z| \rightarrow \infty} V(z) > \inf_{z \in
	\mathbb{R}^N}V(z)=\gamma >0. \eqno{(V_0)}
$$
 In fact, these solutions concentrate at global minimum points of $V$ as  $\epsilon$ tends to 0, c.f. 
 Wang in \cite{W}.  del Pino and Felmer in  \cite{Pino}, assuming a local condition, namely,  there is an 
open and bounded set $\Lambda$ compactly contained  in $\mathbb{R}^{N}$ satisfying
$$
0< \gamma \leq V_0 =\inf_{z\in\Lambda}V(z)< \min_{z \in
	\partial\Lambda}V(z), \eqno{(V_{1})}
$$
 established  the existence of positive  solutions which concentrate around local minimum of $V,$  by introducing 
a penalization method.

\vspace{0.3 cm}
In the nonlocal case, that is, when $s\in (0,1),$ even in the subcritical case,  there are only few references on
 existence and/or concentration phenomena for fractional nonlinear equation $(P_{\epsilon})$, maybe because
 techniques developed for local case can not be adapted imediately, c.f. \cite{Secchi}. For instance, the truncation argument
has to handle carefully in the present situation.  We would like to cite \cite{Secchi, Shang} for the  existence of positive solution, imposing a global condition on $V.$   In \cite{Davila} is studied the existence and concentration phenomena for potential verifying local condition $(V_1)$, and in \cite{Chen, ShangZhang} a concentration phenomenon is treated  near of non degenerate critical point of $V.$

\vspace{0.3cm}
Motivated by approach used in \cite{Pino}, we will establish existence and concentration of positive
 solution for $(P_\epsilon)$,  by supposing that $V$ satisfies $(V_0)-(V_1)$, with $\gamma =V_0$ and 
$V \in L^{\infty}(\mathbb{R}^{N})$. Related to $f$, we will assume the following conditions:\\

\noindent $(f_1)\quad  \, \displaystyle \lim_{t \to 0}\frac{f(t)}{t}=0$ \\

\noindent There is $q \in (2,2^{*}_{s}),$ where $2^{*}_{s}=\frac{2N}{N-2s},$ $N > 2s,$ $s\in (0,1),$ such that
$$
\leqno{(f_2)} \,\,\,\,\, \displaystyle \limsup_{|t| \to +\infty}\frac{|f(t)|}{|t|^{q-1}}<+\infty.
$$

\noindent There is $\theta >2$ such that
$$
\theta F(t) \leq f(t)t, \quad \forall t >0, \leqno{(f_3)}
$$
where $F(t)=\displaystyle \int_{0}^{t}f(\tau)\,d \tau$. \\

\noindent $(f_4)$ \quad The function $\frac{f(t)}{t}$ is increasing for $t>0$.

\vspace{0.5 cm}

Using the above notation we are able to state our main result 

\begin{thm} \label{T2} Assume $(f_1)-(f_4)$ and $(V_0)-(V_1)$, with $\gamma =V_0$ and $V \in
 L^{\infty}(\mathbb{R}^{N})$. Then,  there is $\epsilon_0>0$,  such that $(P_\epsilon)$ has a positive solution
 $u_\epsilon$ for all $\epsilon \in (0,\epsilon_0)$. Moreover, if $x_\epsilon$ denotes a global maximum point of
 $u_\epsilon$, we have 
$$
\lim_{\epsilon \to 0}V(x_\epsilon)=V_0.
$$
  
\end{thm}

The proof of  Theorem \ref{T2} was inspired from   \cite{Pino}, in the sense that we will modify the nonlinearity in a special way and work with an auxiliary problem. Making some estimate we prove that the solutions obtained for the auxiliary problem are solution of the original problem when $\epsilon$ is sufficiently small. However, we would like to point out that we cannot repeat directly the same arguments found in \cite{Pino}. For example, the arguments used in \cite{Pino} do not work well  to prove the Palais-Smale condition for the energy functional associated with auxiliary problem, and also, to show that the solutions of the auxiliary problem are solutions of the original problem for $\epsilon$ small enough. In the present paper, we overcome these difficulties showing a new argument to prove the Palais-Smale condition for the energy functional associated with the auxiliary problem, see proof of Lemma \ref{L2}. Moreover, using some properties of the Bessel kernel, we prove that the obtained solutions of the auxiliary problem are solutions of the original problem for $\epsilon$ small enough, see proofs of Corollary \ref{Cor Limita} and Lemma \ref{Con. ZERO}.

We recall that, for any $s \in (0,1)$, the fractional Sobolev space $H^{s}(\R^N)$ is definied by
\[
H^{s}(\R^N)=\Big\{u\in L^2(\R^N): \ \int_{\R^{2N}}\frac{(u(x)-u(y))^2}{|x-y|^{N+2s}}\ d x\ d y<\infty\Big\},
\]
endowed with the norm
$$
\|u\|_{H^s(\R^N)}=\Big(|u|_{L^2(\mathbb{R}^{N})}^2+\int_{\R^{2N}}\frac{(u(x)-u(y))^2}{|x-y|^{N+2s}}\ d x\ d y\Big)^{1/2}.
$$
The fractional  Laplacian, $(-\Delta)^{s}u,$  of a smooth function $u:\R^{N} \rightarrow  \R$ is defined  by  
$$
{\mathcal F}((-\Delta)^{s}u)(\xi)=|\xi|^{2s}{\mathcal F}(u)(\xi), \ \xi \in \R^N,
$$ 
where ${\mathcal F}$ denotes the Fourier transform, that is, 
\[
{\mathcal F}(\phi)(\xi)=\frac{1}{(2\pi)^{\frac{N}{2}}} \int_{\mathbb{R}^N} \mathit{e}^{-i \xi \cdot x} \phi (x)  \, 
\ d x \equiv \widehat{\phi}(\xi) ,  
\]
for functions $\phi$ in  the  Schwartz class.
Also $(-\Delta)^{s}u$  can be equivalently represented \cite[Lemma 3.2]{nezza} as
$$
(-\Delta)^{s} u(x) = -\frac{1}{2} C(N,s)\int_{\R^N}\frac{(u(x+y)+u(x-y)-2 u(x))}{|y|^{N+2s}}\ d y, \ \forall x \in \R^N,
$$
where 
$$C(N,s)=(\int_{\R^N}\frac{(1- cos\xi_1)}{|\xi|^{N+2s}}d\xi)^{-1},\  \xi=(\xi_1,\xi_2,\ldots,\xi_N).$$
Also, in light of \cite[Propostion~3.4,Propostion~3.6]{nezza}, we have
\begin{equation}
\label{equinorm}
|(-\Delta)^{s/2} u|^2_{L^2(\R^N)}=\int_{\Re^N}|\xi|^{2s}|\widehat{u}|^2d\xi=\frac{1}{2}C(N,s)
\int_{\R^{2N}}\frac{(u(x)-u(y))^2}{|x-y|^{N+2s}}\  d x \ d y,
\end{equation}
for all $u\in H^{s}(\R^N)$,and, sometimes, we identify these two quantities by omitting the normalization constant $\frac{1}{2} C(N,s).$
For $ N > 2s,$  from \cite[Theorem 6.5]{nezza} we also know that, for any $p \in [ 2, 2^{*}_{s}]$,
there exists $C_p>0$ such that
\begin{equation}
\label{emb}
|u|_{L^p(\mathbb{R}^{N})}\leq C_p\|u\|_{H^{s}(\R^N)},
\,\quad\text{for all $u\in H^{s}(\R^N)$}.
\end{equation}

Before to conclude this introduction, we would like point out that using the change variable $v(x)=u(\epsilon x)$, it is possible to prove that $(P_\epsilon)$ is equivalent to the following problem
$$
(-\Delta)^s {u}+V(\epsilon x)u=f(u) \,\,\, \mbox{in} \,\,\, \mathbb{R}^{N}. \eqno{(P_\epsilon)'}
$$

\vspace{0.5 cm}

\noindent \textbf{Notation:} In this paper we use the following
notations:
\begin{itemize}
	\item  The usual norms in $L^{t}(\mathbb{R}^{N})$ and $H^{s}(\mathbb{R}^{N})$ will be denoted by
	$|\,. \,|_{t}$ and $\|\;\;\;\|$ respectively.

	\item   $C$ denotes (possible different) any positive constant.
	
	\item   $B_{R}(z)$ denotes the open ball with center at $z$ and
	radius $R$.
	\item  $a_n \rightharpoonup a$ and $a_n \rightarrow a$  mean the weak and  strong convergence, respectively, as $  n \rightarrow  \infty.$
	
\end{itemize}

\section{The Caffarelli and Silvestre's method}

Hereafter, we will use a method due to Caffarelli and Silvestre in \cite{caffarelli} to get a solution to $(P_\epsilon)'$. In that seminal paper, it was developed a local interpretation of the fractional
Laplacian given in $\R^N$ by considering a Dirichlet to Neumann type operator in the domain
$\R^{N+1}_{+}=\{(x, t) \in \R^{N+1} : t > 0\}.$ A similar extension, in a bounded domain, see for instance, \cite{Brandle, cabre}.
 For $u \in H^{s}(\R^N),$ the solution $w \in X^{1,s}=X^{1,s}(\R^{N+1}_{+})$ of  
 \begin{equation}
 \left\{ \begin{array}{rcl}
 -\dive ( y^{1-2s}\nabla w)=0 & \mbox{in}&  \R^{N+1}_{+}\noindent\\
  w=u& \mbox{on}& \R^{N} \times\{0\}\noindent
 \end{array}\right. 
 \end{equation}
is called $s$-harmonic  extension $w=E_{s}(u)$ of $u$ and it is proved in \cite{caffarelli} that
$$
\lim_{y \to 0^+} y ^{1-2s}\frac{\partial w}{\partial y}(x,y)=-\frac{1}{k_{s}}(-\Delta)^{s}u(x),
$$
where  
$$
k_{s}=2^{1-2s}\Gamma(1-s)/\Gamma (s).
$$ 
Here the spaces $X^{s}(\R^{N+1}_{+})$ and $H^{s}(\R^N)$
are defined as the completion of $C^{\infty}_{0}(\overline{\R^{N+1}_{+}})$ and  $C^{\infty}_{0}(\R^{N}),$
under the norms (which actually do coincide, see \cite[Lemma A.2]{Brandle})
\begin{align*}
\|w\|_{X^{s}}:=&\Big(\int_{ \R^{N+1}_{+}}k_{s}y^{1-2s}|\nabla w|^2 \ud x\ud y\Big)^{1/2},   \\
\|u\|_{H^{s}}:=&\Big(\int_{\R^N}|2\pi\xi|^{2s}|\mathbb{F}(u(\xi))|^2 \ud \xi\Big)^{1/2}=
\Big(\int_{\R^N}|(-\Delta)^{s}u|^2 \ud x\Big)^{1/2}.
\end{align*}
\noindent
Motivated by the above approach  problem, we will study  the existence of solution for the following problem,
 \begin{equation}\label{NPS}
 \left\{ \begin{array}{rcl}
 -\dive ( y^{1-2s}\nabla w)=0 & \mbox{in}&  \R^{N+1}_{+}\noindent\\
  -k_{s}\frac{\partial w}{\partial \nu}=-V(\epsilon x)w + f(w)  & \mbox{on}& \R^{N} \times\{0\},\noindent
 \end{array}\right. 
 \end{equation}
where 
$$
\frac{\partial w}{\partial \nu}=\lim_{y \to 0^+} y ^{1-2s}\frac{\partial w}{\partial y}(x,y),
$$ 
because the function $w \in X^{s}(\R^{N+1}_{+})$ is a solution for this problem, the function $u(x)=w(x,0)$ will be a solution of $(P_\epsilon)'$.

We are looking for a positive solution in the Hilbert space $X^{1,s}$ defined as 
the completion of $C^{\infty}_{0}(\overline{\R^{N+1}_{+}})$ under the norm

$$
\|w\|_{1,s}:= \Big(\int_{ \R^{N+1}_{+}}y^{1-2s}|\nabla w|^2 \ud x\ud y + \int_{\R^N}V(\epsilon x)w(x,0)^2 \ud x\Big)^{1/2}
$$
Here, we are omitting the constant $k_s.$

Associated with (\ref{NPS}), we have the energy functional  $J_\epsilon:X^{1,s} \to \mathbb{R}$ defined by
\begin{equation} \label{E1}
J_\epsilon (v)=\frac{1}{2}\int_{\mathbb{R}_{+}^{N+1}}y^{1-2s}|\nabla v|^{2}\,dxdy+
\frac{1}{2}\int_{\mathbb{R}^{N}}V(\epsilon x)|v(x,0)|^{2}\,dx-\int_{\mathbb{R}^{N}}F(v(x,0))\,dx,
\end{equation}
which is $C^{1}(X^{1,s},\mathbb{R})$ with Gateaux derivative given by

\begin{eqnarray} \label{E2}
J'_\epsilon (v)\phi &=& \frac{1}{2}\int_{\mathbb{R}_{+}^{N+1}}y^{1-2s}\nabla v . \nabla \phi\,dxdy \\ &&+
\frac{1}{2}\int_{\mathbb{R}^{N}}V(\epsilon x) v(x,0) \phi(x,0) \,dx-\int_{\mathbb{R}^{N}}f(v(x,0))\phi (x,0)\,dx \nonumber
\end{eqnarray}
for all $\phi \in X^{1,s}.$

Using some embeddings mentioned in  Br\"andle, Colorado and  S\'anchez \cite{Brandle}, see also \cite{CW,ZhangLiuJiao}, we deduce that the embeddings
$$
X^{1,s} \hookrightarrow L^{p}(\mathbb{R}^{N}) \,\,\, \mbox{for} \,\,\ p \in [2, 2^{*}_{\alpha}] 
$$
are continuous. Moreover, the embeddings 
$$
X^{1,s} \hookrightarrow L^{p}(A) \,\,\, \mbox{for} \,\,\ p \in [2, 2^{*}_{\alpha}) 
$$
are compacts, for any bounded mensurable set $A \subset \mathbb{R}^{N}$.

In what follows,  we will not work directly with functional $J_\epsilon$, because we have some difficulties to prove that it verifies the $(PS)$ condition. Hereafter, we will use the same approach explored in \cite{Pino}, modifying the nonlinearity of a suitable way. 
The idea is the following:  

First of all, without loss of generality, we assume that
$$
V(0)=V_0=\min_{x \in \mathbb{R}^{N}}V(x).
$$
Moreover, since we intend to find positive solutions, we will suppose that
$$
f(t)=0, \quad \forall t \leq 0.
$$
Moreover, let us fix $k>\frac{2 \theta}{\theta -2}$ and $a>0$ verifying
$$
\frac{f(a)}{a}=\frac{V_0}{k},
$$
where $V_0>0$ was given in $(V_0)$. Using these numbers, we set the functions
$$
\tilde{f}(t)=
\left\{
\begin{array}{l}
f(t), \quad t \leq a, \\
\mbox{}\\
\frac{V_0}{k}t, \quad t \geq a
\end{array}
\right.
$$
and
$$
g(x,t)=\chi_\Lambda(x)f(t)+(1-\chi_\Lambda)\tilde{f}(t), \quad \forall (x,t) \in \mathbb{R}^{N} \times \mathbb{R},
$$
where $\Lambda$ was given in $(V_1)$ and $\chi_\Lambda$ denotes the characteristic function associated with $\Lambda$, that is, 
$$
\chi(x)=
\left\{
\begin{array}{l}
1, \quad x \in \Lambda \\
0, \quad x \in \Lambda^c.
\end{array}
\right.
$$
Using the above functions, we will study the existence of positive solution for the following auxiliary problem
$$
\left\{
\begin{array}{l}
(-\Delta)^{s}{u}+V(\epsilon x)u=g_\epsilon( x,u), \quad x \in \mathbb{R}^{N}, \\
\mbox{}\\
u \in H^{1,s}(\mathbb{R}^{N}),
\end{array}
\right.
\leqno{(AP)}
$$
where 
$$
g_\epsilon(x,t)=g(\epsilon x,t), \quad \forall (x,t) \in \mathbb{R}^{N} \times \mathbb{R}.
$$
We recall by using \cite{caffarelli}, the above problem can be  studied as a Neumann problem
$$
 \left\{ \begin{array}{rcl}
 -\dive ( y^{1-2s}\nabla w)=0 & \mbox{in}&  \R^{N+1}_{+}\noindent\\
  -k_{s}\frac{\partial w}{\partial \nu}=-V(\epsilon x)w + g_\epsilon(x, w)  & \mbox{on}& \R^{N} \times\{0\}.\\
 
 \end{array}
\right. 
\leqno{(AP)}
$$
For simplicity we will drop the constant $k_s$ from the above equation.

Here, we would like point out that if $v_\epsilon \in X^{1,s}$ is a solution of $(AP)$ with
$$
v_\epsilon(x,0)<a, \quad \forall x \in \Lambda^{c}_\epsilon,
$$
where $\Lambda_\epsilon=\Lambda / \epsilon$, then  $u_\epsilon(x)=v_\epsilon(x,0)$ is a solution  of $(P_\epsilon)'$.

Associated with $(AP)$, we have the energy functional $E_\epsilon:X^{1,s} \to \mathbb{R}$ given by
$$
E_\epsilon (v)=\frac{1}{2}\int_{\mathbb{R}_{+}^{N+1}}y^{1-2s}|\nabla v|^{2}\,dxdy+
\frac{1}{2}\int_{\mathbb{R}^{N}}V(\epsilon x)|v(x,0)|^{2}\,dx-\int_{\mathbb{R}^{N}}G_\epsilon(x,v(x,0))\,dx
$$
where
$$
G_\epsilon(x,t)=\int_{0}^{t}g_\epsilon(x,\tau)\,d \tau, \quad \forall (x,t) \in \mathbb{R}^{N} \times \mathbb{R}.
$$
Using the definition of $g$, it follows that\\

\noindent $(g_1) \quad \theta G_\epsilon(x,t) \leq g_\epsilon(x,t)t, \quad \forall (x,t) \in \Lambda_\epsilon \times \mathbb{R},$ \\

\noindent and \\

\noindent $(g_2) \quad 2 G_\epsilon(x,t) \leq g_\epsilon(x,t)t \leq \frac{V_0}{k}|t|^{2}, \quad \forall (x,t) \in (\Lambda_\epsilon)^{c} \times \mathbb{R}.$ \\

From $(g_2)$,
$$
L(x,t)=V(x)-G_{\epsilon}(x,t) \geq \left( 1- \frac{1}{2k} \right)V(x)|t|^{2}\geq 0, \quad \forall (x,t) \in (\Lambda_\epsilon)^{c} \times \mathbb{R},
\leqno{(g_3)}
$$
and
$$
M(x,t)=V(x)-g_{\epsilon}(x,t)t \geq \left( 1- \frac{1}{k} \right)V(x)|t|^{2} \geq 0, \quad \forall (x,t) \in (\Lambda_\epsilon)^{c} \times \mathbb{R}.
\leqno{(g_4)}
$$

\begin{lem} \label{L1} The functional $E_\epsilon$ verifies the mountain pass geometry, that is, \\
\noindent $i)$ \quad There are $r,\rho>0$ such that
$$
E_\epsilon(v) \geq \rho, \quad \mbox{for} \quad \|v\|_{1,s}=r
$$
\noindent $ii)$ \quad There is $e \in X^{1,s}$ with $\|e\|_{1,s}>r$ and $E_\epsilon(e)<0$. 
\end{lem}
\noindent {\bf Proof.} \,  From $(g_1)-(g_4)$, there exist $c_1,c_2>0$ verifying
$$
E_\epsilon(v) \geq c_1\|v\|_{1,s}^{2} -c_2\|v\|_{1,s}^{q}, \quad \forall v \in X^{1,s}.
$$
As $q>2$, from the above inequality, there are $r, \rho>0$ such that
$$
E_\epsilon(v) \geq \rho, \quad \mbox{for} \quad \|v\|_{1,s}=r,
$$
showing $i)$. To prove $ii)$, fix $\varphi \in X^{1,s}$ with $supp{\varphi} \subset \Lambda_\epsilon \times \mathbb{R}$. Then, for $t>0$
$$
E_\epsilon(t \varphi)=\frac{t^{2}}{2}\|\varphi\|_{1,s}^{2}-\int_{\mathbb{R}^{N}}F(t\varphi(x,0))\,dx.
$$
From $(f_3)$, we know that there are $c_3,c_4 \geq 0$ verifying
$$
F(t) \geq c_1|t|^{\theta}-c_2, \quad \forall t \geq 0.
$$
Using the above inequality, we derive
$$
\lim_{t \to +\infty}E_\epsilon(t \varphi)=-\infty.
$$
Thereby, $ii)$ follows with $e=t \varphi$ and $t$ large enough.

\fim

\begin{lem} \label{L2} The functional $E_\epsilon$ verifies the $(PS)_c$ condition. 

\end{lem}
\noindent {\bf Proof.} \, Let $(v_n) \subset X^{1,s}$ be a $(PS)_c$ sequence for $E_\epsilon$, that is,
$$
E_\epsilon(v_n) \to c \quad \mbox{and} \quad E'_\epsilon(v_n) \to 0.
$$ 
From this, there are $C_1>0$ and  $n_0 \in \mathbb{N}$, such that
$$
E_\epsilon(v_n)-\frac{1}{\theta}E'_\epsilon(v_n)v_n \leq C_1+C_1\|v_n\|_{1,s}, \quad \forall n \geq n_0.
$$
On the other hand, by $(g_1)-(g_4)$, there is $C_2>0$ such that
$$
E_\epsilon(v_n)-\frac{1}{\theta}E'_\epsilon(v_n)v_n \geq C_2\|v_n\|_{1,s}^{2}, \quad \forall n \in \mathbb{N}.
$$
Thus,
$$
C_2\|v_n\|_{1,s}^{2} \leq C_1+C_1\|v_n\|_{1,s}, \quad \forall n \geq n_0,
$$
showing that $(v_n)$ is bounded in $X^{1,s}$.  As $X^{1,s}$ is reflexive, there is a subsequence of $(v_n)$, still denoted by
 itself, and $v \in  X^{1,s}$ such that
$$
v_n \rightharpoonup v \quad \mbox{in} \quad X^{1,s},
$$
$$
v_n \to v \quad \mbox{in} \quad  L_{loc}^{q}(\mathbb{R}^{N}), \quad \forall q \in [1,2^{*}_{s}) 
$$
and
$$
v_n(x,0) \to v(x,0) \quad \mbox{a.e. in} \quad \mathbb{R}^{N}.
$$
See \cite[Thm 7.1]{nezza}, to see the proof of the above limits.

Using the above limits, it is possible to prove that $v$ is a critical point for $E_\epsilon$, that is, 
$$
E'_\epsilon(v)\varphi=0, \quad \forall \varphi \in X^{1,s}.
$$
Considering $\varphi=v$, we have that $E'_\epsilon(v)v=0$, and so,
\begin{eqnarray*}
&&\int_{\mathbb{R}_{+}^{N+1}}y^{1-2s}|\nabla v|^{2}\,dxdy+\int_{\Lambda_\epsilon}V(\epsilon x)|v(x,0)|^{2}\,dx+
\int_{(\Lambda_\epsilon)^{c}}M(x,v(x,0))\,dx  \\ &&=\int_{\Lambda_\epsilon}f(v(x,0))v(x,0)\,dx
\end{eqnarray*}
On the other hand, using the limit $E'_\epsilon(v_n)v_n=o_n(1)$, we derive that
\begin{eqnarray*}
&&\int_{\mathbb{R}_{+}^{N+1}}y^{1-2s}|\nabla v_n|^{2}\,dxdy+\int_{\Lambda_\epsilon}V(\epsilon x)|v_n(x,0)|^{2}\,dx+
\int_{(\Lambda_\epsilon)^{c}}M(x,v_n(x,0))\,dx \\ && =\int_{\Lambda_\epsilon}f(v_n(x,0))v_n(x,0)\,dx+o_n(1).
\end{eqnarray*}
Since $\Lambda_\epsilon$ is bounded, the compactness Sobolev embedding gives
$$
\lim_{n \to +\infty}\int_{\Lambda_\epsilon}f(v_n(x,0))v_n(x,0)\,dx=\int_{\Lambda_\epsilon}f(v(x,0))v(x,0)\,dx
$$
and
\begin{equation} \label{Eq1}
\lim_{n \to +\infty}\int_{\Lambda_\epsilon}V(\epsilon x)|v_n(x,0)|^{2}\,dx=\int_{\Lambda_\epsilon}V(\epsilon x)|v(x,0)|^{2}\,dx.
\end{equation}
Therefore,
\begin{eqnarray*}
&&\limsup_{n \to +\infty}\left( \int_{\mathbb{R}_{+}^{N+1}}y^{1-2s}|\nabla v_n|^{2}\,dxdy+\int_{(\Lambda_\epsilon)^{c}}M(x,v_n(x,0))
\,dx \right)\\ && =\int_{\mathbb{R}_{+}^{N+1}}y^{1-2s}|\nabla v|^{2}\,dxdy+\int_{(\Lambda_\epsilon)^{c}}M(x,v(x,0))\,dx
\end{eqnarray*}
Now, recalling that $M(x,t) \geq 0$, the Fatous' lemma leads to
\begin{eqnarray*}
&&\
\liminf_{n \to +\infty}\left( \int_{\mathbb{R}_{+}^{N+1}}y^{1-2s}|\nabla v_n|^{2}\,dxdy+\int_{(\Lambda_\epsilon)^{c}}M(x,v_n(x,0))
\,dx \right)\\&& \geq \int_{\mathbb{R}_{+}^{N+1}}y^{1-2s}|\nabla v|^{2}\,dxdy+\int_{(\Lambda_\epsilon)^{c}}M(x,v(x,0))\,dx
\end{eqnarray*}
Hence,
\begin{equation} \label{Eq2}
\lim_{n \to +\infty}\int_{\mathbb{R}_{+}^{N+1}}y^{1-2s}|\nabla v_n|^{2}\,dxdy=\int_{\mathbb{R}_{+}^{N+1}}y^{1-2s}|\nabla v|^{2}\,dxdy
\end{equation}
and
$$
\lim_{n \to +\infty}\int_{(\Lambda_\epsilon)^{c}}M(x,v(x,0))\,dx=\int_{(\Lambda_\epsilon)^{c}}M(x,v(x,0))\,dx
$$
The last limit combined with definition of function $M$ gives
$$
\lim_{n \to +\infty}\int_{(\Lambda_\epsilon)^{c}}V(\epsilon x)|v_n(x,0)|^{2}\,dx=\int_{(\Lambda_\epsilon)^{c}}
V(\epsilon x)|v(x,0)|^{2}\,dx.
$$
Gathering this limit with (\ref{Eq1}), we deduce that
\begin{equation} \label{Eq3}
\lim_{n \to +\infty}\int_{\mathbb{R}^{N}}V(\epsilon x)|v_n(x,0)|^{2}\,dx=\int_{\mathbb{R}^{N}}V(\epsilon x)|v(x,0)|^{2}\,dx.
\end{equation} 
From (\ref{Eq2})-(\ref{Eq3}), we infer that
$$
\lim_{n \to +\infty}\|v_n\|_{1,s}^{2}=\|v\|_{1,s}^{2}.
$$
As $X^{1,s}$ is a Hilbert space and $v_n \rightharpoonup v$ in $X^{1,s}$, the above limit yields
$$
v_n \to v \quad \mbox{in} \quad  X^{1,s},
$$
showing that $E_\epsilon$ verifies the $(PS)_c$ condition.
\fim

\begin{thm} \label{T0} The functional $E_\epsilon$ has a nonnegative critical point $v_\epsilon \in X^{1,s}$ such 
\begin{equation} \label{vepsilon}
E_\epsilon(v_\epsilon)=c_\epsilon \quad \mbox{and} \quad E'_\epsilon(v_\epsilon)=0,
\end{equation}
where $c_\epsilon$ denotes the mountain pass level associated with $E_\epsilon$. 
\end{thm}
\noindent {\bf Proof.} \, The existence of the critical point $v_\epsilon $ is an immediate consequence of
 the Mountain Pass Theorem due to Ambrosetti and Rabinowitz (see e.g.  \cite[Thm 1.17]{willem}). The function $v_\epsilon$ is nonnegative, because
$$
E'_\epsilon(v_\epsilon)(v_{\epsilon}^{-})=0  \Longrightarrow v_{\epsilon}^{-}=0,
$$
where $v_{\epsilon}^{-}=\min\{v_{\epsilon},0\}$.

\fim

\vspace{0.5 cm}

Related to the case $\epsilon=0$, it is possible to prove that there is $w \in X^{1,s}$ such that 
\begin{equation} \label{exist. w}
J_0(w)=c_0 \quad \mbox{and} \quad J'_0(w)=0.
\end{equation}

\vspace{0.5 cm}

The next lemma establishes an important relation between $c_\epsilon$ and $c_0$.

\begin{lem} The numbers $c_0$ and $c_\epsilon$ verify the equality below
$$
\lim_{\epsilon \to 0}c_\epsilon=c_0.
$$ 
\end{lem}
\noindent {\bf Proof.} \, From $(V_0)$,
$$
c_\epsilon  \geq c_0, \quad \forall \epsilon \geq 0.
$$
Then,
\begin{equation} \label{limita1}
\liminf_{\epsilon \to 0}c_\epsilon  \geq c_0.
\end{equation}
Next, fix  $t_\epsilon >0$ such that 
$$
t_\epsilon w \in  \mathcal{M}_\epsilon=\{v \in X^{1,s}\setminus \{0\}\,:\, E'_\epsilon(v)v=0\}.
$$
By definition of $c_\epsilon$, we know that
$$
c_\epsilon \leq \max_{t \geq 0}E_\epsilon(tw)=E_\epsilon(t_\epsilon w).
$$
Now, standard arguments as those used in \cite{AF2}, it is possible to prove that
$$
\lim_{\epsilon \to 0}t_\epsilon=1
$$
and
$$
\lim_{\epsilon \to 0}E_\epsilon(t_\epsilon w)=J_0(w).
$$
Thus, 
\begin{equation} \label{limita2}
\limsup_{\epsilon \to 0}c_\epsilon \leq J_0(w)=c_0.
\end{equation}
From (\ref{limita1})-(\ref{limita2}),
$$
\limsup_{\epsilon \to 0}c_\epsilon =c_0.
$$
\fim

\begin{lem} \label{sequencia} There are $r,\beta, \epsilon^{*}>0$ and $(y_\epsilon) \subset \mathbb{R}^{N}$ such that
$$
\int_{B_r(y_\epsilon)}|v_\epsilon(x,0)|^{2}\,dx\geq \beta, \quad \forall \epsilon \in (0, \epsilon^{*}). 
$$
\end{lem}
\noindent {\bf Proof.} \,  First of all, we recall that since $(v_\epsilon)$ satisfies $(\ref{vepsilon}),$
 there is $\alpha >0$, which is independent of $\epsilon$, such that
\begin{equation} \label{alfa}
\|v_\epsilon\|_{1,s}^{2} \geq \alpha, \quad \forall \epsilon >0.
\end{equation}

To show the lemma, it is enough to see that for any sequence $(\epsilon_n) \subset (0,+\infty)$ with $\epsilon_n \to 0$, the limit below
$$
\lim_{n \to +\infty}\sup_{y \in \mathbb{R}^{N}}\int_{B_r(y)}|v_{\epsilon_n}(x,0)|^{2}\,dx=0,
$$
does not hold for any $r>0$. Otherwise, if it holds for some $r>0$, we can use a Lion's type results for $X^{1,s}$, see \cite{Secchi}, to conclude that
$$
v_{\epsilon_n}(.,0) \to 0 \quad \mbox{in} \quad L^{q}(\mathbb{R}^{N}), \quad \forall q \in (2,2^{*}_{s}).
$$

Using the above limits, it is possible to show that 
$$
\|v_{\epsilon_n}\|^{2}_{1,s} \to 0 \quad \mbox{as} \quad n \to +\infty,
$$
which contradicts (\ref{alfa}). 

\fim

\begin{lem} \label{convergencia em Lambda} For any $\epsilon_n \to 0$, consider the sequence $(y_{\epsilon_n}) 
\subset \mathbb{R}^{N}$given in Lemma \ref{sequencia} and $\psi_n(x,y)=v_{\epsilon_n}(x+y_{\epsilon_n},y).$ Then there is a subsequence of $(\psi_n)$, still denoted by itself, and $\psi \in X^{1,s} \setminus \{0\}$ such that 
\begin{equation} \label{psi}
\psi_n \to \psi \quad \mbox{in} \quad  X^{1,s}.
\end{equation}
Moreover, there is $x_0 \in \Lambda$ such that
\begin{equation} \label{yn}
\lim_{n \to 0}\epsilon_n y_{\epsilon_n}=x_0 \quad \mbox{and} \quad V(x_0)=V_0.
\end{equation} 
\end{lem}
\noindent {\bf Proof.} \, We begin the proof showing that  $(\epsilon_n y_{\epsilon_n})$ is a bounded sequence. Hereafter, we denote by $(y_n)$ and $(v_n)$ the sequences $(y_{\epsilon_n})$ and $(v_{\epsilon_n})$ respectively. 

Since $E'_{\epsilon_n}(v_n)\phi=0, \forall \phi \in X^{1,s}$, we have that
$$
\int_{\mathbb{R}_{+}^{N+1}}y^{1-2s}\nabla v_n \nabla \phi \,dxdy+\int_{\mathbb{R}^{N}}V(\epsilon_n x)v_n(x,0) \phi(x,0)\,dx-
\int_{\mathbb{R}^{N}}g_\epsilon(x,v_n(x,0))\phi(x,0)\,dx=0.
$$
Then, 
$$
\int_{\mathbb{R}_{+}^{N+1}}y^{1-2s}|\nabla v_n|^{2}\,dxdy+\int_{\mathbb{R}^{N}}V(\epsilon_n x)|v_n(x,0)|^{2}\,dx-
\int_{\mathbb{R}^{N}}g_\epsilon(x,v_n(x,0))v_n(x,0)\,dx=0.
$$
From definition of $g$, we know that
$$
g_\epsilon(x,t)\leq f(t), \quad \forall t \geq 0,
$$
and reminding that $v_n \geq 0,$ we infer that 
$$
\int_{\mathbb{R}_{+}^{N+1}}y^{1-2s}|\nabla v_n|^{2}\,dxdy+\int_{\mathbb{R}^{N}}V_0|v_n(x,0)|^{2}\,dx-
\int_{\mathbb{R}^{N}}f(v_n(x,0))v_n(x,0)\,dx\leq0.
$$
Therefore, there is, $s_n \in (0,1)$ such that
$$
s_nv_n \in \mathcal{M}_0=\{v \in X^{1,s}\setminus \{0\}\,:\, J'_0(v)v=0\}.
$$
Using the characterization of $c_0$, we know that
$$
c_0 \leq J_0(s_nv_n), \quad \forall n \in \mathbb{N}.
$$
As 
$$
J_0(w) \leq E_\epsilon(w), \quad \forall w \in X^{1,s} \quad \mbox{and} \quad \epsilon >0,
$$
it follows that
$$
c_0 \leq J_0(s_nv_n) \leq E_{\epsilon_n}(s_nv_n)\leq \max_{s \geq 0}E_{\epsilon_n}(sv_n)=E_{\epsilon_n}(v_n)=c_{\epsilon_n}. 
$$
Recalling that
$$
c_{\epsilon_n} \to c_0,
$$
the last inequality gives
$$
(s_nv_n) \subset \mathcal{M}_0 \quad \mbox{and} \quad J_0(s_nv_n) \to c_0.
$$
By change variable, we also have
$$
(s_n\psi_n) \subset \mathcal{M}_0 \quad \mbox{and} \quad J_0(s_n\psi_n) \to c_0.
$$
Using Ekeland Variational Principle, we can assume that $(s_nv_n)$ is a $(PS)_{c_0}$ sequence, that is,
$$
(s_n\psi_n) \subset \mathcal{M}_0,  \quad J_0(s_n\psi_n) \to c_0 \quad \mbox{and} \quad J'_0(s_n\psi_n) \to 0.
$$
A direct computation shows that $(s_n)$ is a bounded sequence with
$$
\liminf_{n \to +\infty}s_n>0.
$$
Thus, in what follows, we can assume that for some subsequence, there is $s_0>0$ such that
$$
s_n \to s_0.
$$
From definition of $y_n$ and $\psi_n$, we know that $\psi \in X^{1,s} \setminus \{0\}$. Moreover, as $J'_0(s_n\psi_n) \to 0$, 
we also have  $ J'_0(s_0\psi)=0$. Thereby, by definition of $c_0$, 
$$
c_0 \leq J_0(s_0\psi).
$$ 
On the other hand, by Fatous' Lemma
$$
\liminf_{n \to +\infty}J_0(s_n\psi_n) \geq J_0(s_0\psi),
$$
implying that
$$
J_0(s_0\psi)=c_0 \quad \mbox{and} \quad J'_0(s_0\psi)=0.
$$
The above equalities combined with Fatous' Lemma, up to a subsequence,  gives 
$$
s_n \psi_n \to s_0 \psi \quad \mbox{in} \quad X^{1,s}.
$$
Recalling that $s_n \to s_0>0$, we can conclude that
$$
\psi_n \to \psi \quad \mbox{in} \quad X^{1,s},
$$ 
showing (\ref{psi}). 

Using the last limit, we are able to prove (\ref{yn}). To do this, we begin making the following claim

\begin{claim} \label{convergencia yn} $\displaystyle \lim_{n \to +\infty}dist(\epsilon_n y_n, \overline{\Lambda})=0$

\end{claim}

Indeed, if the claim does not hold, there is $\delta>0$ and a subsequence of $(\epsilon_n y_n)$, still denoted by itself, such that,
$$
dist(\epsilon_n y_n, \overline{\Lambda}) \geq \delta, \quad \forall n \in \mathbb{N}.
$$
Consequently, there is $r>0$ such that
$$
B_r(\epsilon_n y_n) \subset \Lambda^{c}, \quad \forall n \in \mathbb{N}.
$$ 
From definition of $\psi_n$, we have that 
\begin{eqnarray*}
&&\int_{\mathbb{R}_{+}^{N+1}}y^{1-2s}|\nabla \psi_n|^{2}\,dxdy+\int_{\mathbb{R}^{N}}V(\epsilon_n x +
\epsilon_n y_n)|\psi_n(x,0)|^{2}\,dx \\ && =\int_{\mathbb{R}^{N}}g(\epsilon_nx +\epsilon_n y_n,\psi_n(x,0))\psi_n(x,0)\,dx.
\end{eqnarray*}
Note that
\begin{eqnarray*}
&&
\int_{\mathbb{R}^{N}}g(\epsilon_nx +\epsilon_n y_n,\psi_n(x,0))\psi_n(x,0)\,dx  \leq 
\int_{B_{\frac{r}{\epsilon_n}}(0)}g(\epsilon_nx +\epsilon_n y_n,\psi_n(x,0))\psi_n(x,0)\,dx \\ && +
\int_{\mathbb{R}^{N} \setminus B_{\frac{r}{\epsilon_n}}(0)}g(\epsilon_nx +\epsilon_n y_n,\psi_n(x,0))\psi_n(x,0)\,dx
\end{eqnarray*}
and so,
\begin{eqnarray*}
&&
\int_{\mathbb{R}^{N}}g(\epsilon_nx +\epsilon_n y_n,\psi_n(x,0))\psi_n(x,0)\,dx \\ &&  \leq \frac{V_0}{k}
\int_{B_{\frac{r}{\epsilon_n}}(0)}|\psi_n(x,0)|^{2}\,dx+\int_{\mathbb{R}^{N} \setminus B_{\frac{r}{\epsilon_n}}(0)}
f(\psi_n(x,0))\psi_n(x,0)\,dx.
\end{eqnarray*}
Therefore,
\begin{eqnarray*}
&&
\int_{\mathbb{R}_{+}^{N+1}}y^{1-2s}|\nabla \psi_n|^{2}\,dxdy+\int_{\mathbb{R}^{N}}V(\epsilon_n x +\epsilon_n y_n)|\psi_n(x,0)|^{2}
\,dx\\ &&  \leq 
\frac{V_0}{k}\int_{B_{\frac{r}{\epsilon_n}}(0)}|\psi_n(x,0)|^{2}\,dx+\int_{\mathbb{R}^{N} \setminus 
B_{\frac{r}{\epsilon_n}}(0)}f(\psi_n(x,0))\psi_n(x,0)\,dx.
\end{eqnarray*}
implying that
\begin{equation} \label{Eq4}
\int_{\mathbb{R}_{+}^{N+1}}y^{1-2s}|\nabla \psi_n|^{2}\,dxdy+A \int_{\mathbb{R}^{N}}|\psi_n(x,0)|^{2}\,dx \leq
\int_{\mathbb{R}^{N} \setminus B_{\frac{r}{\epsilon_n}}(0)}f(\psi_n(x,0))\psi_n(x,0)\,dx,
\end{equation}
where $A=V_0\left(1-\frac{1}{k}\right)$. By (\ref{psi}), 
$$
\int_{\mathbb{R}^{N} \setminus B_{\frac{r}{\epsilon_n}}(0)}f(\psi_n(x,0))\psi_n(x,0)\,dx \to 0
$$
and
\begin{eqnarray*}
&&
\int_{\mathbb{R}_{+}^{N+1}}y^{1-2s}|\nabla \psi_n|^{2}\,dxdy+A \int_{\mathbb{R}^{N}}|\psi_n(x,0)|^{2}\,dx \\ &&  \to 
\int_{\mathbb{R}_{+}^{N+1}}y^{1-2s}|\nabla \psi|^{2}\,dxdy+A\int_{\mathbb{R}^{N}}|\psi(x,0)|^{2}\,dx>0,
\end{eqnarray*}
which contradicts (\ref{Eq4}). This proves Claim \ref{convergencia yn}. 

From Claim \ref{convergencia yn}, there are a subsequence of $(\epsilon_n y_n)$ and $x_0 \in \overline{\Lambda}$ such that
$$
\lim_{n \to +\infty}\epsilon_n y_n=x_0. 
$$
\begin{claim} \label{convergencia yn 2}  $x_0 \in \Lambda$. 
\end{claim}

Indeed, from definition of $\psi_n$,
$$
\int_{\mathbb{R}_{+}^{N+1}}y^{1-2s}|\nabla \psi_n|^{2}\,dxdy+
\int_{\mathbb{R}^{N}}V(\epsilon_n x +\epsilon_n y_n)|\psi_n(x,0)|^{2}\,dx\leq \int_{\mathbb{R}^{N}}f(\psi_n(x,0))\psi_n(x,0)\,dx.
$$
Then, by (\ref{psi}),
$$
\int_{\mathbb{R}_{+}^{N+1}}y^{1-2s}|\nabla \psi|^{2}\,dxdy+\int_{\mathbb{R}^{N}}V(x_0)|\psi(x,0)|^{2}\,dx\leq 
\int_{\mathbb{R}^{N}}f(\psi(x,0))\psi(x,0)\,dx.
$$
Hence, there is $s_1 \in (0,1)$ such that
$$
s_1 \psi \in \mathcal{M}_{V(x_0)}=\{v \in X^{1,s} \setminus \{0\}\,:\,\tilde{J}'_{V(x_0)}(v)v=0 \}
$$
where $\tilde{J}_{V(x_0)}:X^{1,s} \to \mathbb{R}$ is given by 
$$
\tilde{J}_{V(x_0)} (v)=\frac{1}{2}\int_{\mathbb{R}_{+}^{N+1}}y^{1-2s}|\nabla v|^{2}\,dxdy+
\frac{1}{2}\int_{\mathbb{R}^{N}}V(x_0)|v(x,0)|^{2}\,dx-\int_{\mathbb{R}^{N}}F(v(x,0))\,dx.
$$
If $\tilde{c}_{V(x_0)}$ denotes the mountain pass level associated with $\tilde{J}_{V(x_0)}$, we must have 
$$
\tilde{c}_{V(x_0)} \leq \tilde{J}_{V(x_0)}(s_1\psi)\leq \liminf_{n \to +\infty}E_{\epsilon_n}(v_n)=\liminf_{n \to +\infty}c_{\epsilon_n}=c_0=\tilde{c}_{V(0)}.
$$
Hence,
$$
\tilde{c}_{V(x_0)} \leq \tilde{c}_{V(0)}, 
$$
from where it follows that
$$
V(x_0) \leq V(0).
$$
As $V(0)=V_0=\inf_{x \in \mathbb{R}^{N}}V(x)$, the above inequality implies that
$$
V(x_0)=V_0.
$$
Moreover, by $(V_1)$, $x_0 \notin \partial \Lambda$. Then, $x_0 \in \Lambda$, finishing the proof.
\fim

\begin{cor} \label{Cor Limita} Let $(\psi_n)$ the sequence given in Lemma \ref{convergencia em Lambda}. Then, $\psi_n(\cdot,0) \in L^{\infty}(\mathbb{R}^{N})$ and there is $K>0$ such that
\begin{equation} \label{Linfinito}
|\psi_n(\cdot,0)|_{\infty} \leq K, \quad \forall n \in \mathbb{N}
\end{equation}
and
\begin{equation} \label{Linfinito2}
\psi_n(\cdot,0) \to \psi(\cdot,0) \quad \mbox{in} \quad L^{p}(\mathbb{R}^{N}), \quad \forall p \in (2,+\infty).
\end{equation}
As an immediate consequence, the sequence $h_n(x)= g(\epsilon_n x+\epsilon_n y_n,\psi_n(x,0))$ must verify
\begin{equation} \label{Linfinito3}
h_n \to f(\psi(\cdot,0)) \quad \mbox{in} \quad L^{p}(\mathbb{R}^{N}), \quad \forall p \in (2,+\infty).
\end{equation}
\end{cor}

\noindent {\bf Proof. }\, In what follows,  for each $L>0$, we set  
$$
\psi_{n,L}(x,y)=
\left\{
\begin{array}{lcr} 
\psi_n(x,y), &\mbox{if}& \psi_n(x,y)\leq L\\
L, &\mbox{if}& \psi_n(x,y) \geq L
\end{array}
\right.
$$
and
$$
z_{n,L}=\psi_{n,L}^{2(\beta -1)}\psi_n,
$$
with $\beta >1$ to be determined later. Since 
\begin{eqnarray*} &&
\int_{\mathbb{R}_{+}^{N+1}}y^{1-2s}\nabla \psi_n \nabla \phi \,dxdy+\int_{\mathbb{R}^{N}}V(\epsilon_n x +
 \epsilon_n y_n)\psi_n(x,0) \phi(x,0)\,dx \\ && -\int_{\mathbb{R}^{N}}g(\epsilon_n x + \epsilon_n y_n,\psi_n(x,0))\phi(x,0)\,dx=0, \quad \forall \phi \in X^{1,s},
\end{eqnarray*}
adapting the same approach explored in Alves and Figueiredo \cite[Lemma 4.1]{AF2}, we will find the following estimate
$$
|\psi_n(.,0)|_{\infty} \leq C |\psi_n(.,0)|_{{2^{*}_s}}.
$$
As $(\psi_n)$ is bounded in $X^{1,s}$, we conclude that there is $K>0$ such that
$$
|\psi_n(.,0)|_{\infty} \leq K, \quad \forall n \in \mathbb{N}.
$$

Now, the limit (\ref{Linfinito2}) is obtained  by interpolation on the $L^{p}$ spaces, while that (\ref{Linfinito3}) follows 
combining the growth  condition on  $g$ with (\ref{Linfinito2}). 
\fim

In what follows, we denote by $(w_n) \subset H^{s}(\mathbb{R}^{N})$ the sequence $(\psi_n(\cdot,0))$, that is, 
$$
w_n(x)=\psi_n(x,0), \quad \forall x \in \mathbb{R}^{N}.
$$
Since 
\begin{eqnarray*} &&
\int_{\mathbb{R}_{+}^{N+1}}y^{1-2s}\nabla \psi_n \nabla \phi \,dxdy+\int_{\mathbb{R}^{N}}V(\epsilon_n x +
 \epsilon_n y_n)\psi_n (x,0) \phi(x,0)\,dx \\ && -\int_{\mathbb{R}^{N}}g(\epsilon_n x + \epsilon_n y_n,\psi_n(x,0))\phi(x,0)\,dx=0, \quad \forall \phi \in X^{1,s},
\end{eqnarray*}
we have that $w_n$ is a solution of the problem
$$
(-\Delta)^{s}{w_n}+V(\epsilon_n x+\epsilon_n y_n)w_n=g(\epsilon_n x+\epsilon_n y_n,w_n), \quad \mbox{in} \quad \mathbb{R}^{N},
$$
or equivalently,
\begin{equation} \label{C1}
(-\Delta)^{s}{w_n}+w_n=\chi_n, \quad \mbox{in} \quad \mathbb{R}^{N},
\end{equation}
where
\begin{equation} \label{C3}
\chi_n(x)=w_n(x)+g(\epsilon_n x+\epsilon_n y_n x,w_n(x))-V(\epsilon_n x+\epsilon_n y_n)w_n(x), \quad x \in \mathbb{R}^{N}.
\end{equation}
Denoting $\chi(x)=w(x)+f(w(x))-V(x_0)w(x)$, by Corollary \ref{Cor Limita}, we have that
\begin{equation} \label{C3}
\chi_n \to \chi \quad \mbox{in} \quad L^{p}(\mathbb{R}^{N}), \quad \forall p \in [2,+\infty)
\end{equation}
and there is $k_1>0$,
\begin{equation} \label{C4}
|\chi_n|_{\infty} \leq k_1, \quad \forall n \in \mathbb{N}.
\end{equation}

Using some results found in \cite{FQT}, we know that
$$
w_n(x)=\mathcal{K}\ast \chi_n=\int_{\mathbb{R}^{N}}\mathcal{K}(x-y)\chi_n(y)dy,
$$
where $\mathcal{K}$ is the Bessel kernel, which verifies \\

\noindent $(K_1) \quad \mathcal{K} $ is positive, radially symmetric and smooth in $\mathbb{R}^{N} \setminus \{0\}$, \\

\noindent $(K_2)$ \quad There is $C>0$ such that  
$$
\mathcal{K}(x) \leq \frac{C}{|x|^{N+2s}}, \quad \forall x \in \mathbb{R}^{N} \setminus \{0\} 
$$
\noindent and \\

\noindent $(K_3) \quad \mathcal{K} \in L^{q}(\mathbb{R}^{N}), \quad \forall q \in [1,N/N-2s)$. \\

Using the above informations, we are able to prove the following result

\begin{lem} \label{Con. ZERO}The sequence $(w_n)$ verifies 
$$	
w_n(x) \to 0 \quad \mbox{as} \quad |x| \to +\infty,
$$
uniformly in $n \in \mathbb{N}.$
\end{lem}
\noindent {\bf Proof.} Given $\delta >0$, consider the sets
$$
A_\delta=\{y \in \mathbb{R}^{N}\,:\,  |y-x|\geq 1/\delta  \}
$$
and
$$
B_\delta=\{y \in \mathbb{R}^{N}\,:\,  |y-x|<1/\delta  \}.
$$
Hence, 
$$
0\leq w_n(x) \leq  \int_{\mathbb{R}^{N}}\mathcal{K}(x-y)|\chi_n|(y)dy=\int_{A_\delta}\mathcal{K}(x-y)|\chi_n|(y)dy+\int_{B_\delta}\mathcal{K}(x-y)|\chi_n|(y)dy.
$$
From definition of $A_\delta$ and $(K_2)$, we have that, for all $n \in \N,$
\begin{equation} \label{C5}
\int_{A_\delta}\mathcal{K}(x-y)|\chi_n|(y)dy \leq C\delta^{s}|\chi_n|_\infty\int_{A_\delta}\frac{dy}{|x-y|^{N+s}} 
\leq C\delta^{s}k_1\int_{A_1}\frac{dy}{|x-y|^{N+s}}=C_1\delta^{s}.
\end{equation}
On the other hand,
$$
\int_{B_\delta}\mathcal{K}(x-y)|\chi_n|(y)dy\leq \int_{B_\delta}\mathcal{K}(x-y)|\chi_n - \chi|(y)dy+\int_{B_\delta}\mathcal{K}(x-y)|\chi|(y)dy.
$$
Fix $q >1$ with $q \approx 1$ and $q'>2$ such that $\frac{1}{q}+\frac{1}{q'}=1$. From $(K_3)$ and (\ref{C1}), 
$$
\int_{B_\delta}\mathcal{K}(x-y)|\chi_n|(y)dy\leq |\mathcal{K}|_{q}|\chi_n-\chi|_{q'}+|\mathcal{K}|_{q}|\chi|_{L^{q'}(B_\delta)}
$$ 
As 
$$
|\chi_n-\chi|_{q'} \to 0 \quad \mbox{as} \quad n \to +\infty
$$
and
$$
|\chi|_{L^{q'}(B_\delta)} \to 0 \quad \mbox{as} \quad |x| \to +\infty,
$$
we deduce that there are  $R>0$ and $n_0 \in \mathbb{N}$ such that
\begin{equation} \label{C6}
\int_{B_\delta}\mathcal{K}(x-y)|\chi_n|(y)dy\leq \delta, \quad \forall  n \geq n_0 \quad \mbox{and} \quad |x|\geq R.
\end{equation}
From (\ref{C5}) and (\ref{C6}),
\begin{equation} \label{C7}
\int_{\mathbb{R}^{N}}\mathcal{K}(x-y)|\chi_n|(y)dy\leq C_1\delta^{d}+\delta, \quad \forall  n \geq n_0 \quad \mbox{and} \quad |x|\geq R.
\end{equation}
The same approach can be used to prove that for each $n \in \{1,...., n_{0}-1\}$, there is $R_n>0$ such that
\begin{equation} \label{C8}
\int_{\mathbb{R}^{N}}\mathcal{K}(x-y)|\chi_n|(y)dy\leq C_1\delta^{d}+\delta, \quad  |x|\geq R_n.
\end{equation}
Hence, increasing $R, $ if necessary, we must have
$$
\int_{\mathbb{R}^{N}}\mathcal{K}(x-y)|\chi_n|(y)dy\leq C_1\delta^{d}+\delta, \quad \mbox{for} \quad |x|\geq R, \quad \mbox{uniformly in } \quad n \in \mathbb{N}.
$$	
Since $\delta$ is arbitrary, the proof is finished.  \fim

\vspace{0.5 cm}

\begin{cor} \label{Cor original} There is $n_0 \in \mathbb{N}$ such that
$$
v_n(x,0) < a, \quad \forall n \geq n_0 \quad \mbox{and} \quad \forall x \in \Lambda_{\epsilon_n}^{c}.
$$	
Hence, $u_n(x)=v_n(x,0)$ is a solution of $(P_{\epsilon_n})'$ for $n \geq n_0$.
\end{cor}
\noindent {\bf Proof.}\, By Lemma (\ref{convergencia em Lambda}), we know that $\epsilon_n y_n \to x_0 
$, for some $x_0 \in \Lambda$. Thereby, there is $r>0$ such that some subsequence, still denoted by itself,
$$
B_{r}(\epsilon_n y_n) \subset \Lambda, \quad \forall n \in \mathbb{N}.
$$ 
Hence, 
$$
B_{\frac{r}{\epsilon_n}}(y_n) \subset \Lambda_{\epsilon_n}, \quad \forall n \in \mathbb{N},
$$
or equivalently
$$
\Lambda_{\epsilon_n}^{c} \subset B_{\frac{r}{\epsilon_n}}^{c}(y_n), \quad \forall n \in \mathbb{N}.
$$
Now, by Lemma \ref{Con. ZERO}, there is $R>0$ such that
$$
w_n(x) < a, \quad \mbox{for} \quad |x| \geq R \quad \mbox{and} \quad \forall n \in \mathbb{N},
$$
from where it follows,
$$
v_n(x,0)=\psi_n(x-y_n,0)=w_n(x-y_n)<a, \quad \mbox{for} \quad x \in (B_R(y_n))^{c} \quad \mbox{and} \quad \forall n \in \mathbb{N}.
$$
On the other hand, there is $n_0 \in \mathbb{N}$, such that
$$
\Lambda^{c}_{\epsilon_n} \subset (B_{\frac{r}{\epsilon_n}}(y_n))^{c} \subset (B_{R}(y_n))^{c}, \quad \forall n \in \mathbb{N},
$$
then
$$
v_n(x,0)<a, \quad \forall x \in \Lambda^{c}_{\epsilon_n}  \quad \mbox{and} \quad n \geq n_0,
$$
finishing the proof. 
\fim

\section{Proof of Theorem \ref{T2}} 

By Theorem \ref{T0}, we know that problem $(AP)$ has a nonnegative solution $v_\epsilon$ for all $\epsilon>0$.   Applying Corollary \ref{Cor original}, there is $\epsilon_0$ such that
$$
v_{\epsilon}(x,0)<a, \quad \forall x \in \Lambda^{c}_{\epsilon} \quad \mbox{and} \quad \forall \epsilon \in (0, \epsilon_0),
$$
that is, $v_{\epsilon}(\cdot,0)$ is a solution of $(P_\epsilon)'$ for $\epsilon  \in (0, \epsilon_0)$. Considering
$$
u_\epsilon(x)=v_{\epsilon}({x}/{\epsilon},0), \quad \mbox{for} \quad \forall \epsilon  \in (0, \epsilon_0),
$$
it follows that $u_\epsilon$ is a solution for original problem  $(P_\epsilon)$

If $x_\epsilon$ denotes a global maximum point of $u_\epsilon$, it is easy to see that there is $\tau_0>0$ such that 
$$
u_\epsilon(x_\epsilon) \geq \tau_0, \quad \forall \epsilon >0.
$$
In what follows, setting  $z_\epsilon=\frac{x_\epsilon -\epsilon y_\epsilon}{\epsilon}$, we have that $z_\epsilon$ is a global maximum point of $w_\epsilon$ and  
$$
w_\epsilon(z_\epsilon) \geq \tau_0, \quad \forall \epsilon >0.
$$

Now, we claim that
\begin{equation} \label{DV}
\lim_{\epsilon \to 0}V(x_\epsilon)=V_0.
\end{equation}
Indeed, if the  above limit does not hold, there is $\epsilon_n \to 0$ and $\gamma >0$ such that 
\begin{equation} \label{DV0}
V(x_\epsilon)\geq V_0 +\gamma, \quad \forall n \in \mathbb{N}. 
\end{equation}
By Lemma \ref{Con. ZERO}, we know that
$$
w_{\epsilon_n}(x) \to 0 \quad \mbox{as} \quad |x| \to +\infty \quad \mbox{uniformly in} \quad n \in \mathbb{N}.
$$
Therefore, $(z_\epsilon)$ is a bounded sequence. Moreover, for some subsequence, we also know that there is $x_0 \in \Lambda$ satisfying $V(x_0)=V_0$ and 
$$
\epsilon_n y_{\epsilon_n } \to x_0.
$$
Hence,
$$
x_{\epsilon_n}=\epsilon_n z_{\epsilon_n}+\epsilon_n y_{\epsilon_n } \to x_0
$$
implying that
$$
V(x_{\epsilon_n}) \to V_0, 
$$
which is a contradiction with (\ref{DV0}), showing that (\ref{DV}) holds.

\end{document}